\documentclass[10pt]{article}

\usepackage{amsmath}
\usepackage{graphicx}
\usepackage{amsfonts}
\usepackage{amsopn}
\usepackage{amsthm}
\usepackage{verbatim}
\usepackage{ amssymb }

\newcommand{\dm}{d\sigma}

\newcommand{\ba}{{\mathbf{a}}}
\newcommand{\bp}{{\mathbf{p}}}
\newcommand{\bx}{{\mathbf{x}}}
\newcommand{\bm}{{\mathbf{m}}}
\newcommand{\m}{{\mathbf{m}}}

\newcommand{\zt}{\tilde{z}}
\newcommand{\rt}{\hat{\rho}}

\newcommand{\Db}{{\bar{D}}}
\renewcommand{\o}{d}

\newcommand{\Adm}{{\rm{Adm}(\Db,\Tb)}}

\renewcommand{\r}{\rho}

\newcommand{\Cp}{C_{>0}(\bar{D})}

\newcommand{\Ob}{\bar{\Omega}}

\newcommand{\Tb}{{\bar{T}}}

\newcommand{\ra}{\rangle}
\newcommand{\la}{\langle}

\renewcommand{\b}{\beta}
\newcommand{\z}{\zeta}

\newcommand{\F}{\mathcal{F}}

\renewcommand{\epsilon}{\varepsilon}

\renewcommand{\(}{\left(}
\renewcommand{\)}{\right)}

\newcommand{\0}{(0)}
\newcommand{\1}{(1)}
\newcommand{\2}{(2)}
\newcommand{\3}{(3)}

\newtheorem{theorem}{\bf \normalsize \bf  Theorem}[section]

\newtheorem{problem}{\bf \normalsize \bf  Problem}
\newtheorem*{problem'}{\bf \normalsize \bf  Problem I$_{\text{\bf w}}$}

\theoremstyle{remark}

\theoremstyle{definition}
\newtheorem{definition}[theorem]{\bf \normalsize \bf Definition}
\newtheorem{scheme}[theorem]{\bf \normalsize \bf Numerical Scheme}

\title{Iterative scheme for solving optimal transportation problems
arising in reflector design}
\author{Tilmann Glimm and Nick Henscheid\\
Department of Mathematics\\
Western Washington University, Bellingham, WA 98225\\
}
\date{}
\begin{document}
\maketitle
\begin{abstract}
We consider the geometrical optics problem of finding a system of two reflectors that transform a spherical wavefront into a beam of parallel rays
with prescribed intensity distribution.
Using techniques from optimal transportation theory, it has been shown before that this problem is equivalent to an 
infinite dimensional linear programming (LP) problem.
We investigate techniques for constructing the two reflectors numerically. A straightforward discretization of this problem has the disadvantage that the
number of constraints increases rapidly with the mesh size. So with this technique only very coarse meshes are practical.
To address this well-known issue we propose an iterative solution scheme. In each step an LP problem is solved. Information
from the previous iteration step is used to reduce the number of constraints necessary. As a proof of concept we apply our 
proposed scheme to solve a problem with synthetic data. We give evidence that the scheme converges. We also show that it allows for much finer meshes 
than a simple discretization scheme. There exists a growing literature for the application of optimal transportation theory to
other beam shaping problems. Our proposed scheme is easy to adapt for these problems as well.
\end{abstract}
\section{Introduction}
\begin{figure}
\begin{center}
\includegraphics[width=11cm]{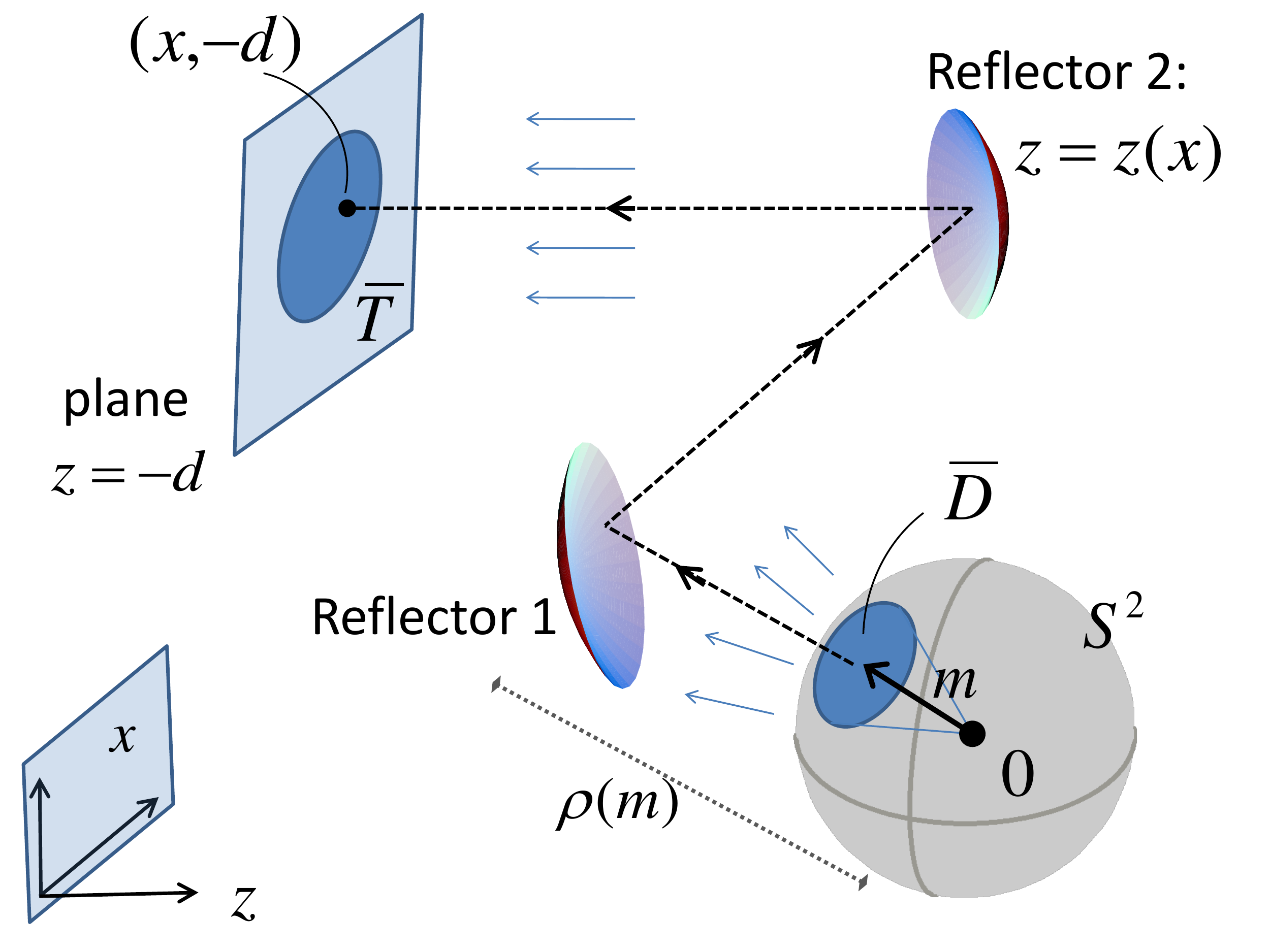}
\end{center}
\caption{Sketch of the reflector problem. The point source is located at the origin $\mathbf{0}$ of the coordinate system.
Note the conventions on coordinates as illustrated 
by the coordinate system in the lower left hand corner: The output beam propagates in the direction
of the {\em{negative}} $z-$axis, and points in the plane perpendicular to the $z-$axis are denoted
by the vector $\bx\in\mathbb R^2$. Thus points in three dimensional space are denoted by $(\bx,z)$. 
See 
Section~\ref{reformulation} for more. (Adapted from \cite{invprob}.)}\label{fig1}
\end{figure}

The following beam shaping problem from geometric optics was described in \cite{OlikerSurvey}; see Figure~\ref{fig1}:
Suppose a point source emits a spherical wavefront with
a given intensity distribution. The problem at hand consists of  transforming this
input beam into an output beam of parallel light rays with a prescribed intensity distribution.
This transformation is to be achieved with a system of two reflectors.
The problem has some practical importance in engineering, see further literature cited in \cite{OlikerSurvey}.  

The first rigorous
mathematical solution to the problem was provided in \cite{invprob}, using an approach
based on the theory of 
optimal transportation \cite{Brenier1,Caf_alloc:96, Gangbo/Mccann:95,villani}.
See also the references \cite{go04, go03, wang, graf, olikerlenses} which deal with other
beam shaping problems using related techniques.

The result of \cite{invprob} is summarized in section~\ref{reformulation} below, with 
Theorem~\ref{thm_equiv} stating the main result.
The central feature is that the original reflector design problem is reformulated as an
infinite dimensional constrained optimization problem, namely the problem of
minimizing a certain functional on a function space. It is the dual problem
for the problem of finding a map from the input aperture to the output aperture
which minimizes a certain cost functional. (See Problem~\ref{probl2} 
in section~\ref{reformulation} for the exact
statement.)

This reformulation of the problem is not only of theoretical value for questions of 
existence and uniqueness
of solutions, but
it also translates into a practical method for finding the solution. In fact,
the discretization of the infinite dimensional constrained optimization problem 
is a standard linear programming problem and can be solved numerically. In this paper, we describe a numerical scheme for solving
these problems. An immediate obstacle is that the number of constraints becomes very large even for relatively coarse meshes -- if there are $M$ sample points
in the input aperture (denoted by $\Db$ in Figure~\ref{fig1})  and $N$ points in the output aperture ($\Tb$ in Figure~\ref{fig1}), then the number of constraints is $M\cdot N$. We found that with a
simple straightforward discretization method, standard LP solvers (on PCs with up to 4 GB RAM) could not handle more than approximately 500 sample points on the domains $\Tb$ and $\Db$ -- 
a number which is arguably too small for many applications. 
For this reason, we devised a more elaborate iterative method, where discretized systems are solved on finer and finer grids, 
in each step using information from the previous solution to choose only a subset of all possible constraints. This
slows the growth of the number of constraints. (Details are described in Section~\ref{sectionscheme}.)
With this step-wise mesh refinement scheme, as a proof of concept,  we were able to obtain solutions on meshes with about 4500 points on each aperture using MATLAB and a standard LP solver.
Using compiled computer language and more specialized solvers for optimal transportation problems, we expect that this number could still be substantially improved.

There has been increased interest in numerical methods for optimal transportation problems in the last 10 years. Most work has concentrated on the
Monge-Amp\`ere equation, which arises in the special case of optimal transport in ${\mathbb{R}}^n$ with quadratic cost (also called $L^2$ optimal transport)
\cite{olikerprussner, Benamoubrenier, benamoufroese, froeseoberman, haberetal, zheligovsky}, although some authors have treated costs proportional to the distance ($L^1$ optimal
transport) \cite{barrettprigozhin}. 
These methods 
are  based on fluid mechanical approaches \cite{Benamoubrenier} or various finite difference approaches  \cite{benamoufroese}. They are generally faster and allow for much larger mesh sizes
than methods based on a discretization of the linear programming 
problem, but since they  use the special 
structure of the  Monge-Amp\`ere equation, they are not directly applicable to more general cost functions or more general manifolds. 

One of the few papers that numerically addressed more general situations for optimal transport is \cite{rueschendorfuckelmann}. A similar approach
to the solution was taken there as we do in this paper, namely a discretization of the linear programming problem. (The authors in \cite{rueschendorfuckelmann}
used the primal formulation of the optimal transportation problem, whereas our approach is equivalent to using the dual formulation.) The authors noted that the
number of variables is growing very fast with mesh sizes, making it impossible to solve the problem numerically even for moderate mesh sizes due to memory
limitations. In fact, they noted that the
software they used, the package Soplex \cite{soplex}, was designed to handle up to 2 million variables. This corresponds to mesh sizes of approximately
700 points on the input and output aperture in our problem if one employs a straightforward discretization scheme. 
The authors noted that for better results, one needs carefully designed programs. This is what we supply here for our problem:
As explained above, our ``mesh refinement''
approach addresses this problem and makes it feasible to solve the problem for finer meshes.

We note that our proposed algorithm does not make any {\it{a priori}} symmetry assumptions on the form the reflectors. 
It can also very easily be adapted for the numerical solution of other beam shaping problems for which a formulation using optimal transportation
theory have been found, for example those in \cite{go03,go04, wang96}.

This article is organized as follows: In section~\ref{reformulation}, we recall the reformulation of the reflector construction problem as a linear
programming transportation as given in \cite{invprob} and fix some notation. We then describe a basic discretization scheme for the numerical solution and 
propose an improved ``mesh refinement'' scheme in section~\ref{sectionscheme}. The following section~\ref{tests} is devoted to numerical tests of this scheme. 
We first derive an explicit
analytical solution to be used as synthetic data for our numerical work. Then we compute the solution numerically and analyze the error of approximation. Our scheme requires a sequence $\epsilon^{(1)}, \epsilon^{(2)}, \ldots$
of ``constraint inclusion thresholds.'' We analyze different choices of this sequence. We conclude with a brief summary and propose future work.

\section{Formulation as Linear Programming Problem}\label{reformulation}
We now briefly review the notation for the problem posed in
\cite{OlikerSurvey}, as well as the result of \cite{invprob}, which reformulates the reflector design problem as a linear programming problem.

Consider the configuration show in Figure~\ref{fig1}. A point source at the origin
$O=(0,0,0)$ generates a spherical wave front over a given input aperture $\Db$ contained in the
unit sphere $S^2$. The input beam has a given intensity distribution. By means of two reflectors, this wave front is to be transformed into a beam of parallel rays
propagating in the direction of the negative $z-$axis. This output beam has to have a prescribed intensity 
distribution.
The cross section of the output beam in a plane perpendicular to the direction of
propagation is called the output aperture, and denoted by $\Tb$. (Certain regularity conditions apply to $\Db$
and $\Tb$, see  \cite{invprob} for the technical details.) 

Denote points in space $\mathbb R^3$ by pairs $(\bx,z)$, where $\bx\in \mathbb R^2$ is the position vector 
in a plane perpendicular to the direction of propagation and $z\in\mathbb R$ is the coordinate in the
(negative) direction of propagation. See again Figure~\ref{fig1} for our convention on the
direction of the $z-$axis. Points on the unit sphere $S^2$ will typically be denoted by $\m\in S^2$;
their components are also written as $\m=(\m_x,m_z)$ with $|\m_x|^2+m_z^2=1$.

We fix the output aperture in the plane $z=-d$. We will seek to represent the first reflector
as the graph of its polar radius $\r(\m)$ (for $\m\in\Db)$, and the second reflector as the graph of a 
function $z(\bx)$ (for $z\in\Tb$). 
(See Figure~\ref{fig1}.) That is
\begin{align}
 {\text{Reflector 1: }}R_1=\{\r(\m)\cdot \m\;\bigl|\; \m\in \Db\},\quad\quad
& {\text{Reflector 2: }}R_2=\{(\bx,z(\bx))\;\bigl|\; \bx\in \Tb\}.    \label{reflectorsets}
\end{align}

The geometrical optics approximation is assumed. It follows from general principles
of geometric optics that  all rays will have equal length from $(0,0,0)$
to the plane $z=-d$; this length is called the optical path length and will be denoted by $L$.
We define the {\em{reduced optical path length}} as $\ell=L-d$. 

Oliker \cite{OlikerSurvey} showed that local energy conservation translates into a complicated
partial differential equation of Monge-Amp\`ere type for $\r(\m)$. As noted in \cite{OlikerSurvey},
the resulting equation is quite involved and  a rigorous analysis of this 
equation seems very difficult. (See equation (59) in \cite{OlikerSurvey}.)

To amend this, the problem was reformulated in \cite{invprob} as a linear programming problem, which makes 
a complete
analysis possible, both concerning theoretical results on existence and uniqueness, and gives a method for practical
computations.

For this, the following function $K(\m,\bx)$, called the cost function in analogy with the theory of
optimal transportation, plays an important role:
\[
   K(\m,\bx)=\frac{\ell - \la \m_x,\bx\ra}{2\ell(\ell^2-|\bx|^2)(1+m_z)}-\frac{1}{4\ell^2} \quad{\text{ for }}\m=(\m_x,m_z)\in\Db, \bx\in\Tb.
\]

In further preparation, the following two transformations are needed:
\begin{definition}\label{def_zt}
Let $z=z(\bx)$ be a continuous function defined on $\Tb\subseteq\mathbb R^2$.
Then define the function
\begin{equation}
\zt(\bx)=\frac{1}{2\ell} - \frac{z(\bx)}{\ell^2-|\bx|^2}\quad\quad \text{   for }\bx\in\Tb.
\end{equation}
\end{definition}

\begin{definition}\label{def_rt}
Let $\r=\r(\m)$ be a continuous function defined on $\Db\subseteq S^2$ with $\r>0$.
Then define the function
\begin{equation}
\rt(\m)=-\frac{1}{2\ell} + \frac{1}{2\r(\m)\cdot (m_z+1)}\quad\quad \text{   for }\m\in\Db.
\end{equation}
\end{definition}

With these preparations, in \cite{invprob}, the following notion of a quasi-reflector pair  and its associated ray tracing map is used.
(The term used in \cite{invprob} is ``reflector pair'', but we use ``quasi-reflector pair'' here for clarity.)
\begin{definition}\label{defreflpair}
A pair $(\r,z)\in\Cp\times C(\Tb)$ is called a {\em{quasi-reflector pair}} if
$\rt,\zt>0$ and
\begin{align}
\rt(\m)& = \sup_{\bx\in\Tb}\(\frac{1}{\zt(\bx)}K(\m,\bx)\) \text{ for } \m\in\Db, \label{def_refl_rt}\\
\zt(\bx)& = \sup_{\m\in\Db}\(\frac{1}{\rt(\m)}K(\m,\bx)\) \text{ for } \bx\in\Tb. \label{refl_zt}
\end{align}
\end{definition}

\begin{definition} \label{def_refl_map}
Let $(\r,z)\in\Cp\times C(\Tb)$ be a quasi-reflector pair. Define its {\em{reflector map}},
or {\em{ray tracing map}}, 
as a set-valued map $\gamma\colon \Db\to\{\text{subsets of }\Tb\}$ via
\[
\gamma(\m)=\{\bx\in\Tb\;\bigl|\; \rt(\m)=\frac{1}{\zt(\bx)}K(\m,\bx)\}\quad \text{ for }\m\in\Db.  
\]
 In \cite{invprob}, it is shown
that $\gamma(\m)$ is in fact single-valued for almost all $\m\in\Db$. Therefore $\gamma$ may be regarded as a  transformation
$\gamma\colon\Db\to\Tb$.
\end{definition}

If $(\r,z)$ is a ``quasi-reflector pair'' in the above sense,the following can be shown \cite{invprob}, justifying the
choice of nomenclature:
If physical copies of the corresponding surfaces (\ref{reflectorsets}) are made from a reflective material,
then 
a ray emitted from the origin in the direction $\m\in\Db$ will be reflected to a ray in the negative $z-$direction
labeled by $\bx=\gamma(\m)\in\Tb$.

With this, the reflector problem can be formulated rigorously:

\begin{problem} {\bf{(Reflector Problem)}}\label{probl1}
For given input and output intensities $I(\bm), \bm\in\Db$,  and $L(\bx), \bx\in\Tb$,  respectively,
satisfying global energy conservation
$\int_{\Db}I(\m)\, \dm =\int_\Tb L(\bx)dx$,
find a pair $(\r,z)\in\Cp\times C(\Tb)$ that satisfies the following conditions:
\begin{enumerate}
\item $(\r,z)$ is a quasi-reflector pair in the sense of Definition~\ref{defreflpair}.
\item The ray tracing map $\gamma\colon\Db\to\Tb$ satisfies \label{loc_energy_pres}
\[
\int_{\gamma^{-1}(\tau)}I(\m)\dm=\int_\tau L(\bx)dx
\]
for any Borel set $\tau\subseteq\Tb$.
\end{enumerate}
\end{problem}

Here $d\sigma$ denotes the standard area element on the sphere $S^2$.
The second condition is local energy conservation.

As we indicate below, the Reflector Problem can be reformulated in the following form:

\begin{problem}\label{probl2}
Minimize the functional
\[
\F(r,\z)=\int_\Db r(\m) I(\m)\dm+\int_\Tb \z(\bx)L(\bx)dx
\]
on the set $\Adm=\{(r,\z)\in C(\Db)\times C(\Tb)\;\bigl|\; r(\m)+\z(\bx)\geq \log K(\m,\bx)
 \text{ for all }\m\in\Db, \bx\in\Tb\}$.
\end{problem}
Problem~\ref{probl2} is the dual problem to the problem of finding a map $P$ from $\Db\subseteq S^2$
to $\Tb\subseteq{\mathbb{R}}^2$ which minimizes the transportation functional 
$P\mapsto \int_{\Db}K(\bm,P(\bm))d\sigma$. 
Problem~\ref{probl1} and Problem~\ref{probl2} are equivalent, as expressed in the following theorem.

\begin{theorem}\cite{invprob} \label{thm_equiv}
Let $(\r,z)\in \Cp\times C(\Tb)$ be a quasi-reflector pair. Then $(\log\rt,\log\zt)\in\Adm$.
The following statements are equivalent:
\begin{enumerate}
\item\label{solp1} $(\r,z)$ solves the Reflector Problem~\ref{probl1}.
\item\label{solp2} $(\log\rt,\log\zt)$ solves Problem~\ref{probl2}.
\end{enumerate}
\end{theorem}
Thus solving the reflector construction problem in Problem~\ref{probl1} 
is equivalent to solving the infinite dimensional LP problem in Problem~\ref{probl2}.
Indeed, it can be shown that a solution as in Theorem~\ref{thm_equiv} exists. (See Corollary 6.5 in \cite{invprob}.)
In the remainder of this paper, we will concentrate on numerical solutions for
this problem.

\section{Numerical Schemes}\label{sectionscheme}
In the following, we describe two scheme for solving Problem~\ref{probl2} numerically.
First, let us list explicitly the given data on which the solution depends:
\begin{enumerate}
\item  domains\footnote{The rigorous result in \cite{invprob} made the additional technical assumption  $(0,0,-1),(0,0,1)\notin\Db$. For practical purposes, these assumptions can often be dropped. }
$\Db\subseteq S^2$, $\Tb\subseteq \mathbb R^2$ 
\item nonnegative integrable functions $I(\bm)$, $\bm\in\Db$, and $L(\bx)$, $\bx\in\Tb$, with $\int_\Db Id\sigma=\int_{\Tb}L dx$
\item a reduced optical length $\ell$
\item the value $\r(\bm_1)=\r_1$ for some fixed $\bm_1\in\Db$, where $\r$ is the radial function of the first reflector.\label{m1}
\footnote{Without this constraint, the reflectors are not uniquely
determined. This can be seen in the formulation of Problem~\ref{probl2}. Indeed, any transformation of the form $r\mapsto r+c$, $\z\mapsto z-c$ for constant $c$
will leave the objective functional unchanged.}
\end{enumerate}

The first numerical scheme, 
Scheme~\ref{brute} is a very straightforward discretization of Problem~\ref{probl2} via meshes 
on the domains $\Db$ and $\Tb$. Because of its simplicity and straightforwardness, we can informally call this a 
``brute force'' technique.  The problem with this scheme is that each pair of points from $\Db$ and $\Tb$
gives rise to a constraint, and so the technique requires a lot of memory,  so that only
coarse meshes can be used. (This problem is very well known, see for example \cite{rueschendorfuckelmann}.)  
To address this issue, we propose a more sophisticated scheme, Scheme~\ref{it},
which is an iterative scheme using finer and finer meshes, where in each iteration, information from
the previous solution is utilized to reduce the number of constraints.

The technique depends on choosing a parameter $\epsilon$, the constraint threshold, in each iteration.
We then discuss the problem of choosing an appropriate sequence of such $\epsilon$.

\subsection{``Simple Discretization'' Scheme} The characterization of solutions by Theorem~\ref{thm_equiv} immediately gives a practical solution scheme for solving the
Reflector Problem~\ref{probl1}, namely by a straightforward discretization of the equivalent Problem~\ref{probl2}. (The same algorithm was used in  
\cite{rueschendorfuckelmann}.)
\begin{scheme} {\bf{(``Simple Discretization'')}}\label{brute}
\begin{itemize}
\item Create a mesh in the input domain $\Db$ by choosing 
sets $\o^{\0}_1,\o^{\0}_2,\ldots,\o^{\0}_{M^{\0}}\subseteq\Db$, where the interiors of any $\o^{\0}_i$ and $\o^{\0}_j$ are disjoint
(for $i\neq j$), and $\displaystyle{\Db=\bigcup_{i=1}^{M^{\0}} \o^{\0}_i}$. (For instance, $\o^{\0}_1,\o^{\0}_2,\ldots,\o_{M^{\0}}$ may
be a triangulation of $\Db$.)
\item Similarly, create a mesh in the output domain $\Tb$ by choosing sets $t^{\0}_1,t^{\0}_2,\ldots,t^{\0}_{N^{\0}}\subseteq\Db$, where the interiors of any $t^{\0}_i$ and $t^{\0}_j$ are disjoint
(for $i\neq j$), and $\displaystyle{\Tb=\bigcup_{j=1}^{N^{\0}} t^{\0}_j}$. 
\item Choose sample points $\bm^{\0}_i\in\o^{\0}_i$ (for $i=1,\ldots,M^{\0}$) and $\bx^{\0}_j\in t^{\0}_j$ (for $j=1,\ldots,N^{\0}$). Here we may assume that $\bm^{\0}_1=\bm_1$ is the same point as given in the data 
in \ref{m1} above.
\item Find the solution $\(\{r^{\0}_i\}_{i=1}^{M^{\0}},\{\z^{\0}_j\}_{j=1}^{N^{\0}}\)$ of the following LP problem:
  \begin{align*}
       & \text{Minimize }
\sum_{i=1}^{M^{\0}}  r^{\0}_i\, I\(\bm^{\0}_i\)\, {\rm{area}}\(d^{\0}_i\)+ \sum_{j=1}^{N^{\0}}  \z^{\0}_j\, L\(\bx^{\0}_j\)\,{\rm{area}}\(t^{\0}_j\)\\
& \text{subject to } r^{\0}_1=\log\rt_1\; (\text{where $\rt_1$ is as given in the data in \ref{m1} above}),\\
& \text{\quad and } r^{\0}_i+\z^{\0}_j\geq \log K\(\bm^{\0}_i,\bx^{\0}_j\) \text{ for }i=1,\ldots, M^{\0};\; j=1,\ldots, N^{\0}.
  \end{align*}           
\item Find the numbers $\rho^{\0}_i, i=1,\dots, M^{\0}$ and $z^{\0}_j, j=1,\ldots, N^{\0}$ such that $r^{\0}_i=\log\rt^{\0}_i$, and $\z^{\0}_j=\log\zt^{\0}_j$.
(This is straightforward by taking the inverse of the transformations given in Definitions~\ref{def_zt} and  \ref{def_rt}.)
\end{itemize}
\end{scheme}
Then $\rho^{\0}_i$ is an approximation for the true value $\rho\(\bm^{\0}_i\)$ of the radial function of the first reflector, evaluated at the sample point $\bm^{\0}_i$ for $i=1,\ldots,M^{\0}$.
Similarly, $z^{\0}_j$ is an approximation for the true value $z\(\bx^{\0}_j\)$ of the function describing the second reflector, evaluated  at the sample point $\bx^{\0}_j$ for $j=1,\ldots,N^{\0}$. (We are using the superscript ``(0)'' here to distinguish this solution from additional iterative approximations we will obtain with the 
iterative scheme described below.)

It is important to point out that this discretization also yields a discretized version of the ray tracing map $\gamma$ (see Definition~\ref{def_refl_map}). Namely, for each index $i=1,\ldots, M^{\0}$, there is at least
one corresponding index $j^{*}$, $1\leq j^{*}\leq N^{\0}$ where the constraint corresponding to the pair $(i,j^*)$ is active, that is, where we have
\[
     r^{\0}_i+\z^{\0}_{j^*}= \log K\(\bm^{\0}_i,\bx^{\0}_{j^*}\) .
\]
This means that the point $\bm^{\0}_{j^{*}}$ is approximately the image of the point $\bx_i^{\0}$ under the ray tracing map $\gamma$.

\subsection{``Iterative Refinement'' Scheme} \label{sectionit}
As noted in the introduction, one of the drawbacks of scheme~\ref{brute} is that the constraint set in the linear programming problem
in the penultimate step becomes large very fast with finer meshes, and the corresponding problem becomes too large to handle for standard LP solvers. We were not
able to solve problems with more than very roughly 1000 mesh points on a standard PC with 4 GB RAM ($M^{\0}\approx 500, N^{\0}\approx 500$), which is arguably too small for many applications.

We addressed this problem by developing an iterative scheme. First the problem is solved for a mesh with relativey few sample points (say $M^{\0}=N^{\0}=250$).
Then a finer mesh is chosen, and the previous solution is used to reduce the number of constraints. 

Specifically, we have the following scheme, depending on a number $\epsilon^{\1}>0$, which we call the ``inclusion threshold,'' to be explained 
in detail below.

\begin{scheme} {\bf{(``Iterative Refinement'')}}\label{it}
\begin{itemize}
\item Use scheme~\ref{brute} to find an initial solution  $(\{r^{\0}_i\}_{i=1}^{M^{\0}},\{\z^{\0}_j\}_{j=1}^{N^{\0}})$ for $M^{\0}$ sample points on $\Db$
and $N^{\0}$ sample points on $\Tb$, respectively.
\item Create a mesh on $\Db$ by choosing sets $\o^{\1}_1,\o^{\1}_2,\ldots,\o^{\1}_{M^{\1}}\subseteq\Db$, where the interiors of any $\o^{\1}_i$ and $\o^{\1}_j$ are disjoint
(for $i\neq j$), and $\displaystyle{\Db=\bigcup_{i=1}^{M^{\1}} \o^{\1}_i}$.
Here $M^{\1}$ is chosen larger than $M^{\0}$, meaning that we have a finer mesh.
\item Similarly, create a mesh on $\Tb$ by choosing sets $t^{\1}_1,t^{\1}_2,\ldots,t^{\1}_{N^{\1}}\subseteq\Db$, where the interiors of any $t^{\1}_i$ and $t^{\1}_j$ are disjoint
(for $i\neq j$), and $\displaystyle{\Tb=\bigcup_{j=1}^{N^{\1}} t^{\1}_j}$. 
\item Choose sample points $\bm^{\1}_i\in\o^{\1}_i$ (for $i=1,\ldots,M^{\1}$) and $\bx^{\1}_i\in t^{\1}_j$ (for $j=1,\ldots,N^{\1}$). Here we may assume that $\bm^{\1}_1=\bm_1$ is the same point as given in the data 
in \ref{m1} above.
\item Interpolate the initial solution $\{r^{\0}_i\}_{i=1}^{M^{\0}}$ to find an approximation for the first reflector. That is, by some interpolation
method, find a continuous function $r(\bm), \bm\in\Db$, with
\[
         r(\bm_i^{\0})=r_i^{\0}\quad\quad \text{for }i=1,\ldots,M^{\0}.
\] 
\item Similarly, interpolate the initial solution $\{\z^{\0}_j\}_{j=1}^{N^{\0}}$ to find an approximation for the second reflector. That is, by some interpolation
method, find a continuous function $\z(\bx), \bx\in\Tb$, with
\[
         \z(\bx_j^{\0})=\z_j^{\0}\quad\quad \text{for }j=1,\ldots,N^{\0}.
\] 

\item Find the solution $\(\{r^{\1}_i\}_{i=1}^{M^{\1}},\{\z^{\1}_j\}_{j=1}^{N^{\1}}\)$ of the following LP problem:
  \begin{align*}
       & \text{Minimize }
\sum_{i=1}^{M^{\1}}  r^{\1}_i\, I\(\bm^{\1}_i\)\, {\rm{area}}\(d^{\1}_i\)+ \sum_{j=1}^{N^{\1}}  \z^{\1}_i\, L\(\bx^{\1}_j\)\,{\rm{area}}\(t^{\1}_j\)\\
& \text{subject to } r^{\1}_1=\log\rt_1,\\
& \text{\quad and } r^{\1}_i+\z^{\1}_j\geq \log K\(\bm^{\1}_i,\bx^{\1}_j\) \\
& \quad\quad\quad\quad\text{ for those pairs } (i,j) \text{ with } r(\bm^{\1}_j)+\z(\bx_j^{\1})-\log K\(\bm^{\1}_{i},\bx^{\1}_{j}\)< \epsilon^{\1}.
  \end{align*}      

\item Find the numbers $\rho^{\1}_i, i=1,\dots, M^{\1}$, and $z^{\1}_j, j=1,\ldots, N^{\1}$, 
such that $r^{\1}_i=\log\rt^{\1}_i$, and $\z^{\1}_j=\log\zt^{\1}_j$.
\end{itemize}
\end{scheme}

The idea of the second scheme  \ref{it} is to solve the discretized LP on a coarse mesh first and then use this information
to reduce the number of constraints needed for the LP on a finer mesh.
To wit, the difference between schemes~\ref{brute} and \ref{it} is that in the discretized LP problem, in scheme~\ref{it} not {\em{all}} pairs 
of sample points from $\Db$ and $\Tb$, represented by pairs of indices  $(i,j)$, are included in the list of constraints.
In fact, we would in principle only need to include those where the constraint is active, that is, where the corresponding inequality holds with equality.
Of course, this information is not available {\em{a priori}}. Instead, we use the following heuristic:
A constraint should only be included in the LP problem if it is ``almost'' active
when we use an interpolation of the solution on the coarse mesh as an approximate solution on the fine mesh. 
This is where the ``inclusion threshold'' $\epsilon^{\1}$ come into play. Namely, we only include those 
constraints $(i,j)$ for which the difference  $ r(\bm^{\1}_j)+\z(\bx_j^{\1})-\log K\(\bm^{\1}_{i},\bx^{\1}_{j}\)$ is less than the inclusion threshold $\epsilon^{\1}$.
(Here $r(\bm)$ and $\z(\bx)$ are interpolations of the solution of the LP problem for the coarse mesh.)
Clearly, increasing $\epsilon^{\1}$ means  more constraints are included, but also potentially the solution of the LP represents a better approximation
of the exact reflector pair.

In our numerical tests, we found it advantageous to scale the functions $I(\bm)$ and $L(\bx)$ so that the approximations
for the integrals $\int_\Db Id\sigma$ and $\int_Tb L dx$ yield exactly the same result.

Note that the algorithm in scheme~\ref{it} can be iterated again.   We can use the first iterative solution $\rho^{\1}_i, i=1,\dots, M^{\1}$, and 
$z^{\1}_j, j=1,\ldots, N^{\1}$, in Scheme~\ref{it} to obtain a second iterative solution $\rho^{\2}_i, i=1,\dots, M^{\2}$, and 
$z^{\2}_j, j=1,\ldots, N^{\2}$, and so on.

\subsection{Choice of thresholds for Iterative Refinement scheme~\ref{it}} \label{sectionthresholds}
For the iteration of scheme~\ref{it} described above, we have to choose a corresponding sequence of mesh sizes $\{(M^{\0}, N^{\0}), (M^{\1}, N^{\1}), (M^{\2}, N^{\2}), \ldots\}$.
(Here $M^{(k)}$ denotes the number of sample points for the input aperture $\Db$, and $N^{(k)}$ denotes the number of sample points for the target aperture $\Tb$, both in the 
$k^{th}$ iteration of the scheme.) We also need to choose a corresponding sequence of thresholds $\epsilon^{\1}, \epsilon^{\2}, \epsilon^{\3}, \ldots$ .

Let us assume that the sequence $\{(M^{\0}, N^{\0}), (M^{\1}, N^{\1}), (M^{\2}, N^{\2}), \ldots\}$ is given. This raises the question: How should the
corresponding sequence of thresholds $\epsilon^{\1}, \epsilon^{\2}, \epsilon^{\3}, \ldots$ be picked?

This is a question of great practical importance. To wit, if the sequence $\{\epsilon^{(k)}\}_k$ is constant or decreases too slowly, then ``many'' constraints are included in each iteration,
meaning that the LP problems will become large fast and we will run out of memory quickly after a few iterations. 
If the sequence $\{\epsilon^{(k)}\}_k$ decreases too fast, then ``few'' constraints are included in each iteration. This could cause for instance that some index $i$ (with $1\leq i \leq M^{(k)}$) 
may not even be included in any of the constraints with any of the indices $j$ (with $1\leq j \leq N^{(k)})$, or vice versa. This would 
cause the LP problem to be unbounded and hence there would be no solution. (This behavior is illustrated in the practical tests summarized in Table~\ref{tableCa}; see also 
the discussion in section~\ref{tests}.)


We can make a very rough estimate for a good choice of the sequence $\{\epsilon^{(k)}\}_k$, assuming for simplicity that $M^{(k)}\approx N^{(k)}$.
Let us assume that each sample point $m_i^{(k)}$ in the input aperture $D$ is paired via the ray tracing map $\gamma$ with a unique point $\bx_{j(i)}^{(k)}\approx \gamma(m_i^{(k)})$ and vice versa.
Thus there are approximately $M^{(k)}$ such pairs.
Let $r(\bm), \bm\in\Db$, and $\z(\bx), \bx\in\Tb$, denote the function obtained by interpolating the solution of the $(k-1)^{th}$ iterate, 
$\(\{r^{(k-1)}_i\}_{i=1}^{M^{(k-1)}},\{\z^{(k-1)}_j\}_{j=1}^{N^{(k-1)}}\)$.
Then the set of all points in $(\bm,\bx)-$space that satisfy $r(\bx)-z(\bx)-\log K(\bm,\bx)<\epsilon^{(k)}$ is approximately an ellipsoid for small $\epsilon^{(k)}$. 
This can be seen by using the Taylor expansion around the pair $(\bm_i^{(k)},\bx_{j(i)}^{(k)})$ in $(\bm,\bx)-$space.
 The volume of this ellipsoid is proportional
to $\left(\epsilon^{(k)}\right)^2$. Since each pair $(\bm_i^{(k)},\bx_j^{(k)})$ that satisfies this inequality corresponds to a constraint included in the LP of the
$k^{th}$ iteration, the total number of constraints should be roughly proportional to $M^{(k)}\cdot(\epsilon^{(k)})^2$. So to keep the number of constraints approximately constant
throughout the iterations, based on these heuristics, one may  choose
\[
   \epsilon^{(k)}=\frac{C}{\sqrt{M^{(k)}}},
\]
where $C$ is a constant. These heuristics are very rough and in practice the number of constraints is still 
increasing with $M^{(k)}$, albeit at a much slower rate than the case of constant $\epsilon^{(k)}$. In section~\ref{results}, we describe an implementation of this idea.

\section{Numerical Tests}\label{tests}
To test the validity of the numerical scheme described above, we used it on a case where the solution is known in analytic form.
In the next subsection, we first describe this special analytic solution, and then discuss our results.
\subsection{An Analytic Solution}
In order to obtain a data set where an explicit analytic solution is known, we
consider a special given configuration of reflectors, and then solve the ``forward'' problem of determining the 
output intensity this system produces for a given input intensity.

\subsubsection{Construction of a pair of reflectors $R_1$, $R_2$}
\begin{figure}
  \begin{center}
\includegraphics[width=12cm]{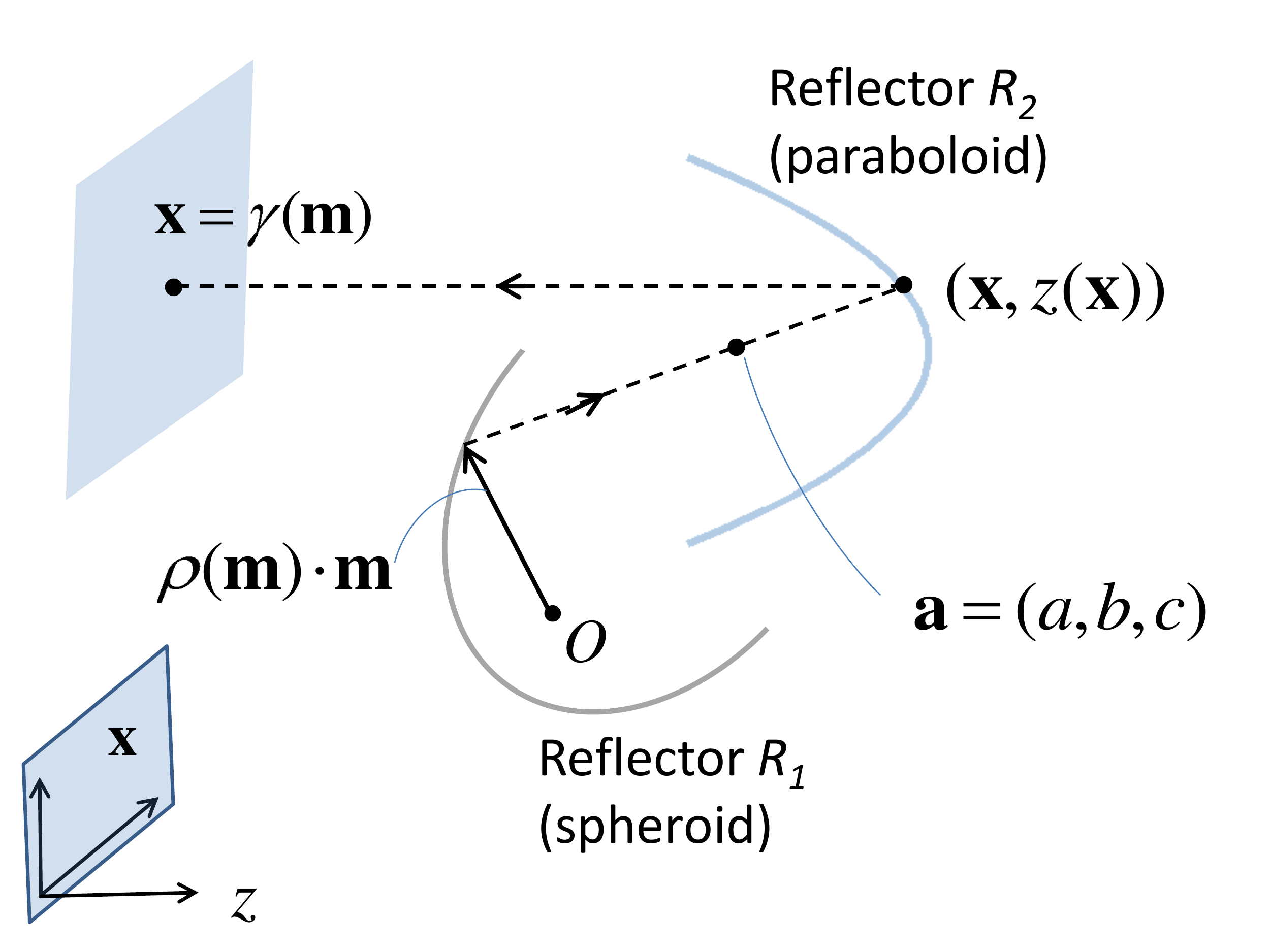}
\end{center}
\caption{Sketch of the reflectors $R_1$ and $R_2$. $R_1$ is the boundary of
an ellipsoid whose foci are the points $O$ and $\ba$. $R_2$ is the boundary of a paraboloid
whose focus is $\ba$.
The sketch shows a two-dimensional cross section. 
A ray given by the direction $\bm\in S^2$
will be reflected by $R_1$ and then by $R_2$ to a ray in the direction of the negative $z-$axis.
}\label{figrefl}
\end{figure}

Consider the following set of two reflectors $R_1$ and $R_2$, sketched in Figure~\ref{figrefl}: Let $\ba=(a,b,c)$ be a given point in $\mathbb{R}^3$,
and let $R>0$ and $\alpha>0$ be two positive numbers with $R>|\ba|=\sqrt{a^2+b^2+c^2}$. 
Now let the first reflector $R_1$ be the boundary of a spheroid
with foci at the origin $\mathbf{O}$ and at $\ba$ and with major diameter $R/2$. (Thus for each point $\mathbf{P}$
on $R_1$, the sum of the distances from each of the two foci $\overline{\mathbf{P}\mathbf{O}}+\overline{\mathbf{P}\ba}$
equals $R$.)
Let the second reflector $R_2$ be the boundary of a paraboloid whose main axis is the negative $z-$axis and whose focus is 
at $\ba$ with focal parameter $2\alpha$.
(Thus for each point $\mathbf{p}=(x,y,z)$ on $R_2$, the sum of the distance to the focus and the shifted $z-$component 
$\overline{\mathbf{P}\ba}+(z-c)$ equals $2\alpha$.)

Note that by definition of the two reflectors, any cone of light rays emitted from the origin $\bf{O}$ will be transformed
into a beam of parallel rays traveling in the direction of the negative $z-$axis. Indeed, a ray emitted from $\mathbf{O}$
will be reflected off $R_1$ towards the focus $\ba$. This ray will then be reflected in the direction of the negative $z-$axis
by $R_2$. See again Figure~\ref{figrefl}.

It is not hard to find explicit expressions for the two reflectors. If we write $R_1=\{ \rho(\bm)\cdot\bm\,\bigl|\, \bm\in S^2\}$
and $R_2=\{ (\bx, z(\bx))\, \bigl|\, \bx=(x,y)\in \mathbb{R}^2\}$, then we have the following expressions:
\begin{align}\label{eqrefl1}
        \text{radial function $\rho$ for $R_1\colon$ }\rho(\bm)=\frac{R^2-|\ba|^2}{2(R-\bm\cdot\ba)}
\end{align}
and
\begin{align}\label{eqrefl2}
        \text{equation for $R_2\colon$ }-4\alpha q=|\bp|^2-4\alpha^2
\end{align}
where we used the shifted coordinates $(\bp,q)=(x,y,z)-(a,b,c)=(\bx, z)-\ba$.

\subsubsection{Ray Tracing Map $\gamma$}
One can now determine the ray tracing map $\gamma$ corresponding
to the reflector pair $R_1$, $R_2$. See again Figure~\ref{figrefl}.

Let $\bm\in S^2$ be a given direction. Thus $\bm=(\bm_\bx,m_z)$ with $\bm_\bx\in\mathbb{R}^2$ and $|\bm_\bx|^2+m_z^2=1$.
To determine the image $\gamma(\bm)$, consider the ray emitted from the origin in the direction $\bm$.
The ray will encounter the first reflector at the point $\rho(\bm)\cdot\bm$ and will then be reflected
towards the second focus $\ba$. Thus the reflected ray can be parameterized by
$     \ba+\lambda\cdot (\ba-\rho(\bm)\cdot\bm)$,
where $\lambda$ is a parameter. 
With equation (\ref{eqrefl2}), this yields that the point where the reflected ray hits the second reflector is $\ba+\lambda^*\cdot (\ba-\rho\cdot\bm)$
where
\[
  \lambda^*=\frac{2\alpha}{(R-\rho)^2-(c-\rho m_z)^2}\left((R-\rho)-(c-\rho m_z)\right)=\frac{2\alpha}{R+c-\rho(1+m_z)}.
\]
The projection of this point to the $xy-$plane is the value of the ray tracing map. Thus we have the following
explicit formula for the ray tracing map:
\begin{align}\label{reflmap}
\gamma(\bm)=\(\begin{array}{c}  a \\b\end{array} \)+ \frac{2\alpha}{R+c-\rho(\bm)\,(1+m_z)} \(\(\begin{array}{c}  a \\b\end{array} \)-\rho(\bm)\, \bm_\bx\),\quad \bm=(\bm_\bx,m_z)\in S^2. 
\end{align}

\subsubsection{Solution of the Forward Problem}
We can now solve the forward problem: Given the reflector pair $R_1, R_2$ as defined above,
an input aperture $\Ob\subseteq S^2$ and an intensity distribution $I(\bm)$, $\bm\in\Ob$, find 
the output aperture $\Tb\subseteq\mathbb{R}^2$ and the output intensity $L(\bx)$, $\bx\in\Tb$.

The output aperture is simply the image of $\Ob$ under the ray tracing map $\gamma$:
$\Tb=\gamma(\Ob)$. To find the induces intensity $L(\bx)$, use the defining property
\begin{align}\label{energybalance}
        \int_\omega I(\bm)d\sigma=\int_{\gamma(\omega)}L(\bx)d\bx
\end{align}
for all Borel sets $\omega\subseteq\Ob$. This is an energy balance equation.

The above integral equation allows us to find an explicit expression for $L(\bx)$ given $I(\bm)$, or vice versa.
We consider for simplicity the case that $\Ob$ is contained in the left hemisphere
$S^2_-=S^2\cap\{(x,y,z)\in\mathbb{R}^3\,\bigl|\, z<0\}$. We can then use coordinates
\begin{align}\label{deftau}
   \tau\colon \{(m_x,m_y)\in\mathbb{R}^2\,\bigl|\, m_x^2+m_y^2<1\}\to S^2_-,\\
               (x,y)\mapsto \(m_x,m_y,-\sqrt{1-m_x^2-m_y^2}\).
\end{align}
In these coordinates, the standard measure on $S^2$ is given by $\displaystyle{d\sigma=\frac{dm_x\,dm_y}{|m_z|}}$,
where\newline $\displaystyle{m_z=-\sqrt{1-m_x^2-m_y^2}}$.
Using this, (\ref{energybalance}) is then equivalent to the equation
\begin{align}\label{IL}
  I(\bm_\bx,m_z)\,\frac{1}{|m_z|}=L(\gamma(\tau(\bm_\bx)))\,J(\gamma\circ\tau).
\end{align}
Here $\displaystyle{J(\gamma\circ\tau)=\left|\det\(\begin{array}{c c} 
       \frac{\partial (\gamma\circ\tau)}{\partial m_x} & \frac{\partial (\gamma\circ\tau)}{\partial m_y} 
               \end{array}\)\right|}  $ denotes the Jacobian of the map $\gamma\circ\tau$,
With the help of a computer algebra system like Mathematica, this can be evaluated explicitly as
\begin{align}\label{Jacobian}
   J(\gamma\circ\tau)=\frac{4\alpha^2(|\ba|^2-R^2)^2}{-m_z \(2(c+R)\,\bm_\bx\cdot\(\begin{array}{c}a \\ b\end{array}\)-(1+m_z)|(a,b)|^2-
                                                     (c+R)^2(1-m_z)\)^2},
\end{align}
where again $\displaystyle{m_z=-\sqrt{1-m_x^2-m_y^2}}$.

\subsection{A data set} \label{dataset}
\begin{figure}
\includegraphics[width=7cm]{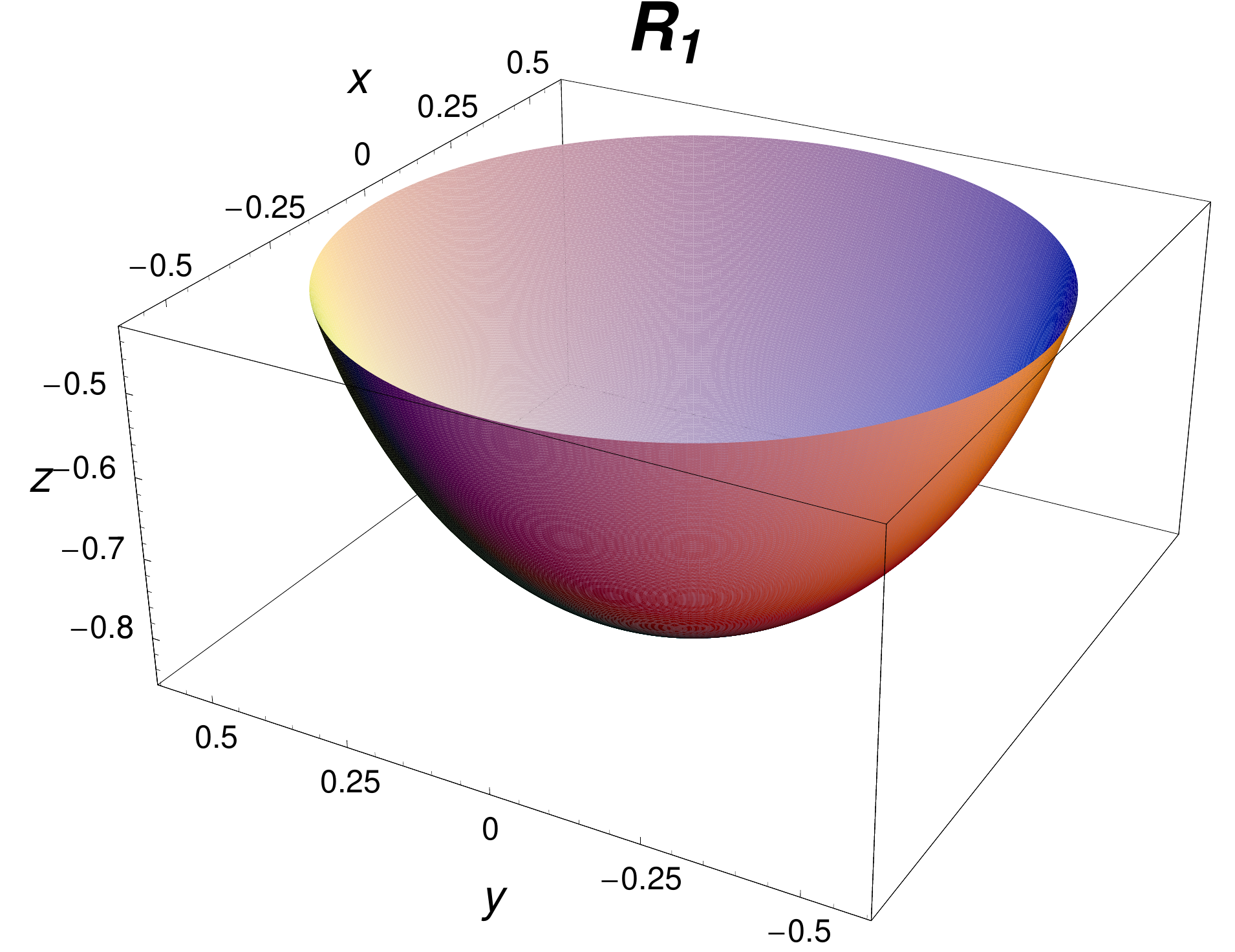} 
\includegraphics[width=7cm]{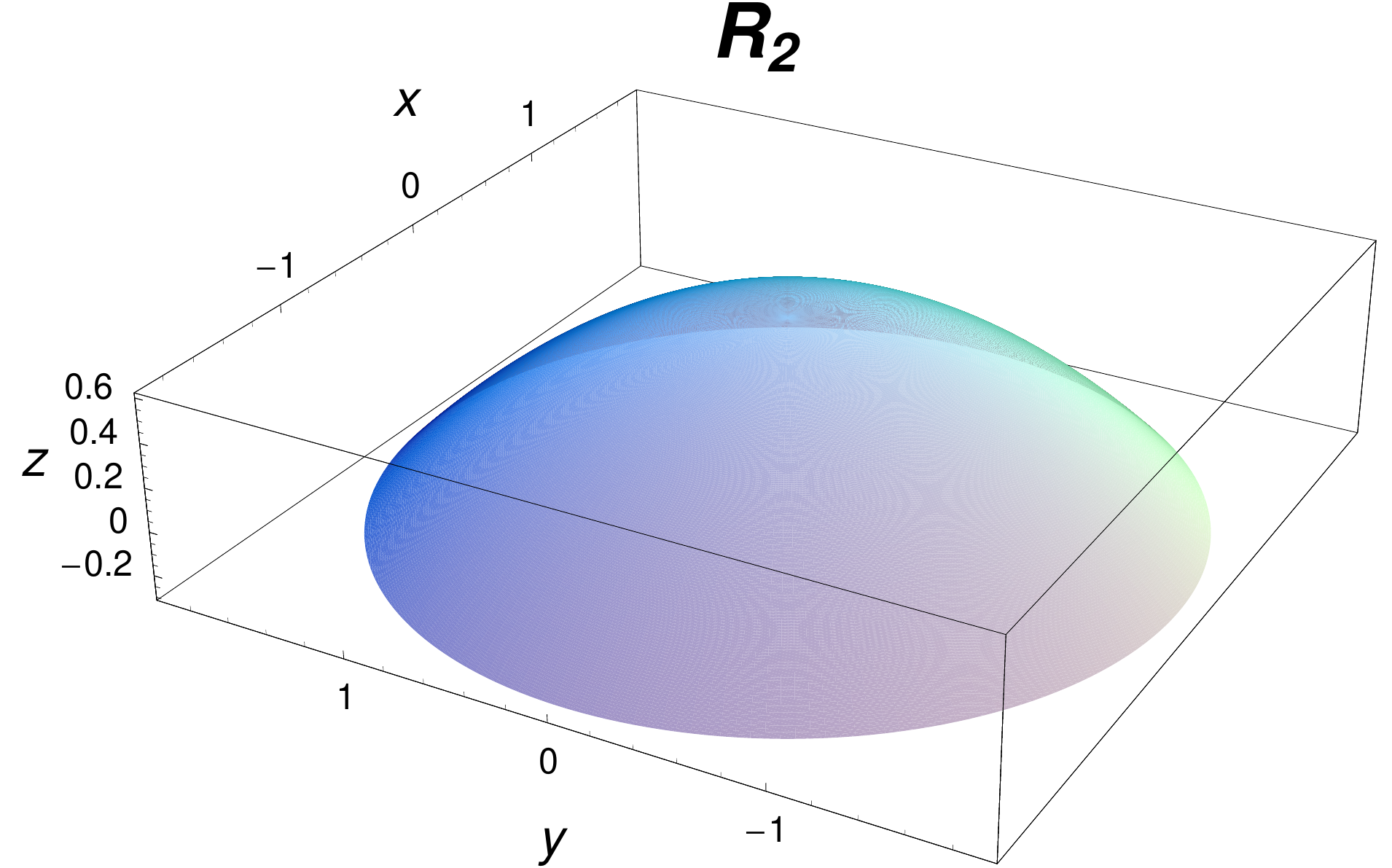} 
\caption{Plots of the reflectors given by (\ref{exact_R1}) (left) and (\ref{exact_R2}) (right) in the explicit data set given in section~\ref{dataset}.
}\label{figanarefl}
\end{figure}

We can now construct a particular data set using the results of the previous sections. We pick $a=b=0$, and
$c=-0.4$, $\alpha=1$, $R=1.3$,
and consider the input aperture 
\begin{align*}
\Ob=\{(m_x,m_y,m_z)\subseteq S^2\,\bigl|\, \sqrt{m_x^2+m_y^2}\leq 0.8, m_z<0 \},
\end{align*}
which is a spherical cap 
centered at the point  $(0,0,-1)$ 
\footnote{
The data do not represent a physically possible set of reflectors since there would be blockage. However, the example is useful as a proof of concept,
and the numerical scheme itself is completely independent of the shape of the apertures or any {\it{a priori}} symmetries of the problem.
}
. Using the results from the previous
subsections, we can now find the output aperture in a straightforward manner:
\begin{align*}
\Tb=\{ (x,y)\in\mathbb{R}^2\,\bigl|\, \sqrt{x^2+y^2}\leq 1.8888888889  \}.
\end{align*}
We choose a constant output intensity
\begin{align*}
L(\bx)&=1\quad \text{ for }\bx\in\Tb.
\end{align*}
Using the relation (\ref{IL}), this gives the input intensity
\begin{align*}
I(\bm)&=\frac{14.2716049383}{\(1-m_z\)^2}\quad \text{ for }\bm=(m_x,m_y,m_z)\in\Ob.
\end{align*}
The two reflectors
are given as follows:
\begin{align}
        \text{radial radius $\rho$ for $R_1\colon$ }\rho(\bm)=\frac{0.765}{1.3+0.4\, m_z}  \label{exact_R1}
\end{align}
and
\begin{align}
        \text{equation for $R_2\colon$ } z=-0.25\,\(x^2+y^2\)+0.6. \label{exact_R2}
\end{align}
The corresponding reduced optical path length is
\[
  \ell=R+2\alpha=2.9.
\]
This follows from the construction of the two reflectors and the properties of ellipsoids and paraboloids.
The two reflectors are shown in Figure~\ref{figanarefl}.

\subsection{Results} \label{results}
Using the data set from subsection~\ref{dataset}, we conducted a series of numerical tests to establish the validity
of the numerical scheme \ref{it}. Since we have explicit analytical expressions for the two solutions
$\rho(m)$ and $z(x)$ describing the two reflectors as given in (\ref{exact_R1}) and (\ref{exact_R2}), respectively,
we were able to compute the error of approximation directly.  (See also Figure~\ref{figanarefl} for plots of these reflectors.)

The iterative scheme \ref{it} was implemented in MATLAB with the package 
distmesh \cite{distmesh} for mesh generation in each iteration step and lp\textunderscore solve \cite{lpsolve} for solving the resulting sparse LP.
Computations were performed on an Intel I5 Dual Core/2.53GHz with 4GB RAM.
Our primary aim is to
demonstrate the practicability of the proposed scheme.
While there are several ways in which the  implementation could be made more efficient, for instance by using a compiled
computer language or a more specialized LP solver for transportation problems, we believe that our results show that the proposed scheme
is easy to implement and makes it possible to use much finer meshes than a simple discretization scheme. 

The mesh generation algorithm employed by the software we used is based on a relaxation scheme of forces in a truss structure \cite{distmesh}.
For this, user input is the desired edge length $h$, not the number of mesh points $M$ (for the domain $\Ob$) or $N$ (for the domain $\Tb$), respectively. 
However, it is not hard to see that the relation between
the average edge length $h$ and the number of mesh points $M$ on $\Db$ obeys
\begin{equation} \label{hN}
   h\sim \frac{1}{\sqrt{M}}.
\end{equation}
Indeed, the total area $A$ is proportional to the number of mesh points $M$ times the area of one triangle of a Delaunay triangulation. The area of a triangle
is proportional to the square of the average edge length $h$. The same relation
holds of course for the number of mesh points $N$ on $\Tb$. (In Tables~\ref{tableCa} and \ref{tablecalibrate}
in this section, we show the number of mesh points $M, N$ instead of
the edge length $h$, as we believe the former are more informative.)  

We chose a sequence of desired edge lengths $h^{\0}, h^{\1},\ldots$. Then, based on the heuristics given above and in section~\ref{sectionthresholds},   
we determined the corresponding  constraint thresholds with the formula

\begin{equation} \label{epsilon}
  \epsilon^{(k)}=C\cdot \left( h^{(k)}\right)^a\quad\quad (k=1,2,3,\ldots),
\end{equation}
where $C$ and $a$ are constants.
We tested several values for the proportionality constant $C$ and the exponent $a$. 
(The heuristic arguments in section~\ref{sectionthresholds} yields $a=1$, but we also tested other values.) Some results are summarized in Table~\ref{tableCa}.
As indicated in section \ref{sectionthresholds}, ``large'' values for $C$  mean that ``many'' 
constraints are included in the LP in each step,
which potentially means good approximations, but also high memory usage, and thus we may run out of memory 
after a few iterations.
In contrast, ``small'' values for $C$ improve memory usage, but may mean that the numerical results may be less 
accurate approximations of the exact results, or indeed the problem may become unbounded, as typically
seems to happen in practice. (See Table~\ref{tableCa}.) A similar logic holds for the exponent $a$:
For larger $a$, the threshold $\epsilon^{(k)}$ decreases faster, leading to less memory usage,
but the danger of the problem becoming unbounded.

Table~\ref{tableCa} confirms and illustrates these principles.
At each iteration, the LP is either unbounded or 
there is a solution. Interestingly, if there is a solution, it seems largely independent of the number of 
constraints, as indicated by the fact that the error of approximation varies very little. (See fourth
column in Table~\ref{tableCa}.) This makes some intuitive sense as the most important issue in each iteration
is really that we  have included the active
constraints. As long as all active constraints are included, we get the same solution. 
If however we did not include
any of the constraints for a certain mesh point, the problem becomes unbounded.

In fact, Table~\ref{tableCa} illustrates the following feature of the algorithm, which is intuitively quite obvious,
although we do not give a formal proof: 
At each iteration step, there is a critical value $\epsilon_{\rm{crit}}$ such that  if 
$\epsilon<\epsilon_{\rm{crit}}$, then the LP is unbounded because no constraints involving a certain point are included.
If  $\epsilon>\epsilon_{\rm{crit}}$, the LP has a solution. As indicated above, it also appears that if there is a solution,
it is likely typically close to the solution of the ``full'' problem (that is, the problem including
all possible constraints). This leads to a possible alternative algorithm for choosing the thresholds
$\epsilon^{(1)}, \epsilon^{(2)}, \ldots$ : Instead of choosing each $\epsilon^{(k)}$ with a pre-determined
formula, use the corresponding critical value $\epsilon_{\rm{crit}}$  in each iteration step. This critical
value can be found by simple trial and error. While this idea is promising, it is possible that such a scheme
may lead to a much slower decay of the maximum error, or even an increase of the error. We did not test this 
idea in this article, but it would be interesting to analyze it further.

\begin{table}
\begin{tabular}{| l | l | l | l | l | l |}
    \hline
$C$   &    $a$ &  sequence of mesh points &    notes  &  max no. constr. & \% of possible constr. \\
\hline
1    &    1.1   &    284;455;724   &   unbounded on 724   &    N/A & N/A \\
1.2  &    1.1   &   284;455;724;1147  &       unbounded on 1147 &    N/A & N/A\\ 
1.3   &   1.1   &   284;455;724;1147 &        Max err 0.0018575  &       222220 & 16.9\% \\
1.5    &   1.1  &    284;455;724;1147    &      Max err 0.0018521   &     258185 & 19.6\%\\
2     &   1.1   &   284;455;724;1147    &     Max err 0.0018521   &     342846 & 26.1\% \\
5    &    1.1   &   284;455;724;1147    &     Max err 0.0018521   &     713571 & 54.0\% \\
0.8  &    1     &   284;455;724  &    unbounded on 724  &    N/A &  N/A\\
0.9  &    1     &   284;455;724;1147    &     unbounded on 1147 &    N/A & N/A\\ 
1    &    1     &   284;455;724;1147   &      Max err 0.0018521  &      226231 &  17.1\% \\
1.5  &    1      &    284;455;724;1147   &       Max err 0.0018521     &    340038 & 25.7\% \\
1    &     1.5    &  284;455      &     unbounded on 455  &    N/A  & N/A \\
2     &   1.5     &  284;455      &  unbounded on 455   &   N/A   & N/A \\
3      &  1.5   &   284;455;724;1147       &  unbounded on 1147  &   N/A  & N/A\\
3.5   &   1.5     &  284;455;724;1147     &    unbounded on 1147   &  N/A   & N/A \\
4      &  1.5  &    284;455;724;1147     &    Max err 0.0018521   &     224140  & 16.9\% \\
5    &    1.5   &   284;455;724;1147     &    Max err 0.0018521     &   282415  & 21.4\% \\
\hline
    \end{tabular}
\caption{Table summarizing the results for runs of the iterative scheme with
different values for the constants $C$ and $a$ given in (\ref{epsilon}). The third column shows the sequence
of mesh points $M^{(k)}$ for $\Db$. (So for instance, in the first row, in the initial discretization, there 
were 284 points. In the first iteration, there were 455 mesh points. The second iteration had 724 points,
but the problem was unbounded.) 
The third column notes if the
problem became unbounded before reaching 1147 mesh points. If the total run was completed up to
1147 mesh points, the third column gives the maximum norm $\| \rho_{\rm{num}}-\rho_{\rm{ana}}\|_{\infty}$ of the error between the resulting numerical solution $\rho_{\rm{num}}$
and the analytic solution $\rho_{\rm{ana}}$. 
The fifth column shows the number of constraints for the case of 1147 mesh points,
if the run reached this maximum number. The last column indicates the corresponding percentage of 
number of constraints relative to the total number of constraints in the ``brute force'' scheme~\ref{brute}.
(So for example in the last row, there were only 21.4\% of all possible constraints used in formulating the last LP,
thus effectively reducing the number of constraints needed by 78.4\%.)
}\label{tableCa}
\end{table}

\begin{table}
\begin{tabular}{| l | l | l | l | l | l | l |}
    \hline
$C$   &    $a$ & no. iterations & no. mesh points  & notes &  no. constr. & \% of possible constr. \\
\hline
0.8  &    1  & 2  &   724    &  Max err 0.00148   &    101174 & 19.4\%\\ 
0.8  &    1  & 3  &   1148    &     Unbounded &  178976   & 13.6\%\\ 
1    &    1  & 3   & 1148  &     Max err 0.0018521  &      226231 &  17.2\% \\
1.2  &    1.1 & 3  &   1148  &       Max err 0.0018575 & 203535   & 15.5\%\\ 
1.3   &   1.1 & 3  &   1148 &        Max err 0.0018575  &       222220 & 16.9\% \\
1.5    &   1.1 & 3  &    1148    &      Max err 0.0018521   &     258185 & 19.6\%\\
2     &   1.1 & 3  &   1148    &     Max err 0.0018521   &     342846 & 26.1\% \\
1    &    1  & 4   &  1824   &     Unbounded  &    446239 &  13.5\% \\
1.5  &    1 &  4   &   1824  &     Max err 0.00131     &    688163 & 20.7\% \\
1.6  &    1 &  4   &   1824  &     Max err 0.00131     &    733894 & 22.05\% \\
1.5  &    1 &  5   &  2882  &     Max err 0.00069      &   1373976 & 16.5\% \\
1.6  &    1  &  5  &   2882&    Max err 0.00069   &  1472275 & 17.7\%\\
1.5  &    1 &  6   &   4536  &       Unbounded    &    2650000 & 12.8\% \\
1.6  &    1  &  6  &   4536&    Max err 0.00067   &  2856627 & 13.88\%\\
1.7  &    1  &  6  &   4536&    Max err 0.00067   &  3057070 & 14.86\%\\
1.6  &    1  &  7  &   7130 &     Out of memory   &  N/A & N/A\\
1.7  &    1  & 7   &   7130 &     Out of memory   &  N/A & N/A\\
\hline
    \end{tabular}
    \caption{Table illustrating the process of calibrating $C$ and $a$ from (\ref{epsilon}) for maximum mesh sizes.  Each row represents a different run of the iterative numerical scheme. 
The third column is the number $k$ of iterations in the run and the fourth is the number $M^{(k)}$ of points in the mesh for $\Db$ in the last iteration.  The fifth column notes the result of the run -- if a solution was found, the maximum error 
is given.  If not, the linear program either became unbounded, or the computer ran out of memory.  The last two columns are the number of constraints used for the last iteration and the percentage of total possible constraints (compared to the ``brute force'' scheme~\ref{brute}).}
\label{tablecalibrate}
\end{table}

To find out the maximum mesh size on which a solution could be obtained without running out of memory or arriving at an unbounded problem,
we ``calibrated'' the values of $C$ and $a$. The process is summarized in Table~\ref{tablecalibrate}.

\begin{figure}
\includegraphics[width=7cm]{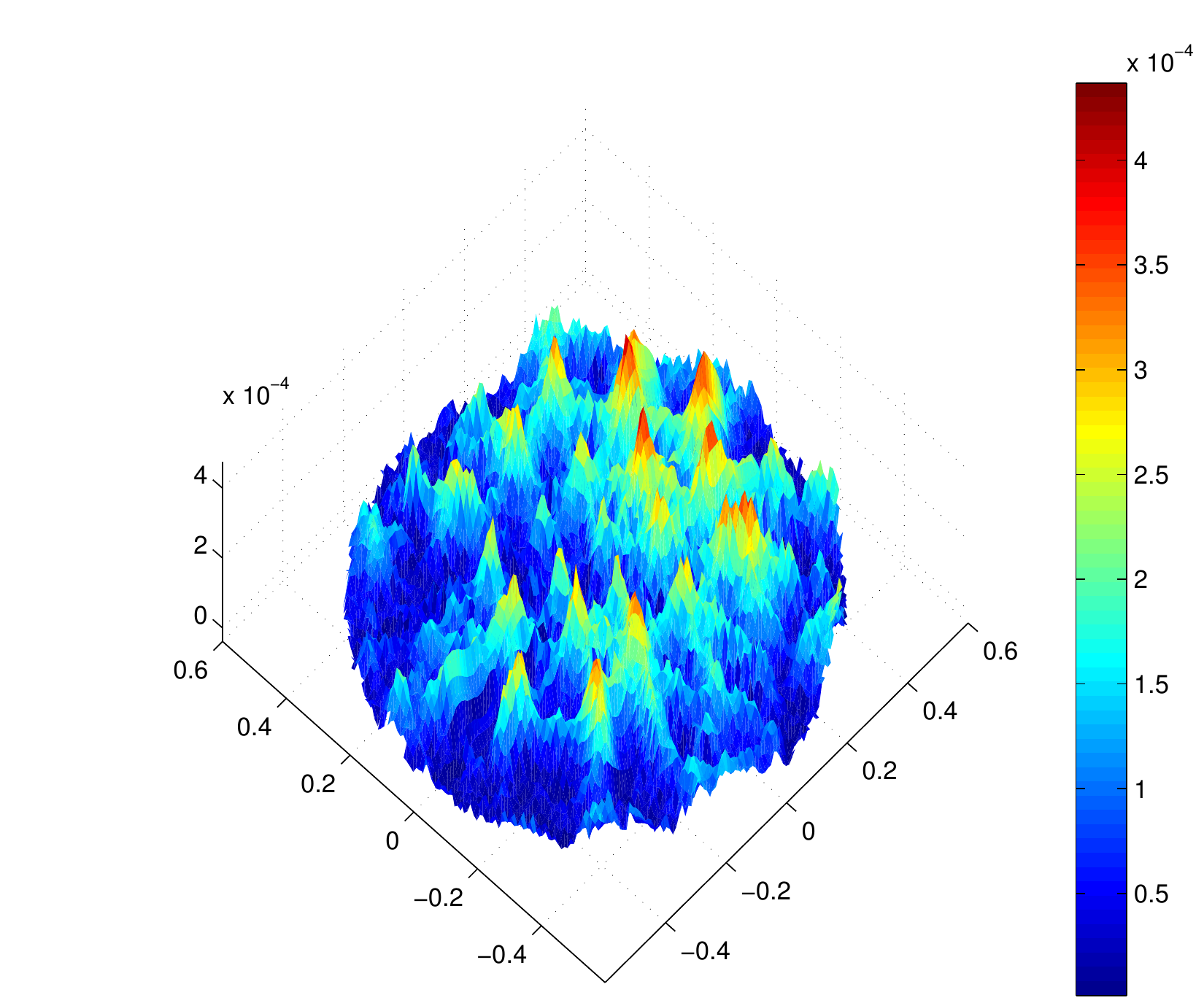} 
\includegraphics[width=7cm]{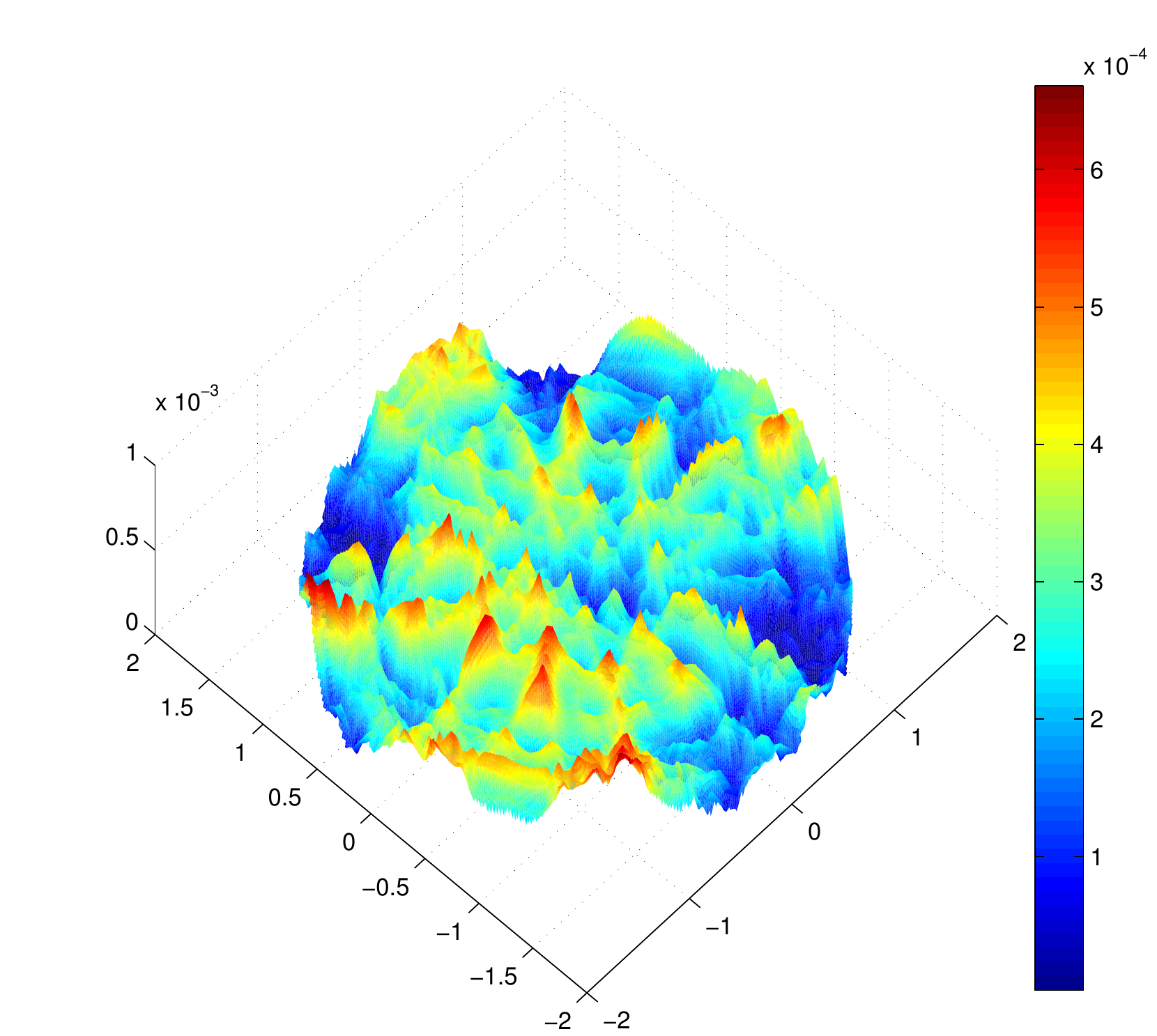} 
\caption{Error surface plots for iteration 7 in Table~\ref{table_error} for reflector 1 (left) and reflector 2 (right). For better visibility, the 
plots were obtained by interpolating the 
error at each mesh point to a continuous function.  Note that the error distribution appears to be relatively uniform, although the ``spikes'' are noticably
arranged in triangular patterns and in a fashion somewhat reminiscent of spokes on a wheel. This may be an artifact of an iterated error in conjuncture with triangular meshes.
}\label{figure_error}
\end{figure}

We obtained best results (that is, largest mesh sizes) for the run with  $C=1.7$ and  $a=1$.
Table~\ref{table_error} now summarizes the results for this specific run in more detail. 
The table shows the errors between the 
numerical computed approximations and the 
exact formulas as given in (\ref{exact_R1}) and (\ref{exact_R2}). Figure~\ref{figure_error} shows 
plots of the approximation errors for the reflectors as  functions on the input and the output apertures.

These results are a good indication that the scheme converges, and that the error decreases
proportional to $M^{-\alpha}$, or equivalently proportional to $h^{2\alpha}$, where roughly 
$0.5\lesssim \alpha\lesssim 1$.
(The results from Figure~\ref{tableCa} appear to indicate $\alpha\simeq 1$, but in other data, we also 
encountered results where $\alpha$ seemed to be closer to $0.5$.) 
The results also illustrate the practicability of the scheme. Note that the number of mesh 
points was increased by a factor of about 16  from iteration 1 to iteration 6. In fact the problem in iteration 6 would
be impossible to solve with a simple discretization scheme due to the size of the problem. 
In iteration 6, we were able to reduce the number of constraints used to only about 15\% of the number of constraints used for the ``naive'' simple discretization scheme~\ref{brute}.

\begin{table}
\begin{tabular}{| l | l | l | l | l | l | l | l | l |}
    \hline
   $k$ &  $h$ &  $M$ & $N$ & Constraints & Max error $R_1$ & $L_2$ error $R_1$ & Max error $R_2$ & $L_2$ error $R_2$ \\ \hline
    1& 0.12 & 284 & 278 & 78,952 & 0.0048 & 0.00143 & 0.008 & 0.0021 \\ \hline
    2& 0.096 & 455 & 450 & 85942 & 0.0022 & 0.00076 & 0.0047  & 0.0014 \\ \hline
    3 & 0.0768 & 724 & 721 & 183308 &0.00148 & 0.00056 & 0.0039 & 0.0012\\\hline
    4 & 0.06144 & 1148 & 1146 & 381971 &0.0012 & 0.00039& 0.00185 &0.00044 \\\hline
    5 & 0.049152 & 1824 & 1810 & 779053 &0.00060 & 0.00021& 0.0013 &0.00033 \\
    \hline
    6& 0.0393216 & 2882 & 2879 & 1,568,533 &0.00059 & 0.00019 & 0.00069 & 0.00016 \\
    \hline
    7 & 0.0314573 & 4536 & 4525 & 3,057,070 &0.00045 & 0.00010 &0.00067 & 0.00027\\
    \hline
    \end{tabular}
\caption{Error for a run of the iterative scheme ~\ref{it} with $C=1.7$ and $a=1$ in (\ref{epsilon}). 
First column: iteration number. Second solumn: average edge length $h$. Third and fourth columns: number of mesh points
for the domains of reflector 1 and reflector 2, respectively. Fifth column: resulting number of constraints in the  LP. 
The scheme stopped after iteration 7 because there was not enough memory available for 
the next iteration. The maximum error as a function of total mesh points $N_{\rm{tot}}=N+M$ decayed in the form $E_{\max}\sim (N_{\rm{tot}})^{\alpha}$ with 
$\alpha\approx -0.82$ for reflector 1 and $\alpha\approx -0.95$ for reflector 2. (These values for $\alpha$ were obtained through least squares curve fitting.)
}\label{table_error}
\end{table}

\section{Conclusion and Future Work}
We proposed two numerical schemes for solving an infinite dimensional optimal transportation problem arising
in reflector design: A straightforward discretization and an improved iterative scheme, which uses knowledge
of the previous solution in each step to reduce the number of constraints. This scheme is easily
adapted to similar transportation problems arising in beam shaping problems, e.g. in \cite{go03} and \cite{go04}. We showed that this 
new scheme is easy to implement and makes it possible to solve the problem on much finer meshes.

There are a number of possible directions for future research. We did not rigorously prove that the scheme converges,
although we strongly expect that it does. In fact the decreasing error shown in Table~\ref{table_error} gives evidence for this, at least in our example.
It would be valuable to have a rigorous proof of convergence.

Another direction for further investigation is to change  the selection algorithm for the threshold
values $\epsilon^{(1)}, \epsilon^{(2)}, \ldots$ as suggested in section~\ref{results}: 
It is a straightforward conjecture that in each step, there is a critical threshold inclusion number
such that the problem is unbounded for values of $\epsilon$ below that number.
One way to select a value of $\epsilon$ is to pick a value for $\epsilon^{(k)}$ close to the  critical value $\epsilon_{\rm{crit}}$.
It would be interesting to flesh out this ideas and analyze the results more closely.

Another avenue of research is to see whether our proposed solution  algorithm can be made parallizable.

It would also be worthwhile to investigate whether the existing PDE-based schemes for  solving the transportation cost with quadratic costs
\cite{Benamoubrenier, benamoufroese, froeseoberman, zheligovsky}
can be adapted to the problem at hand.

\subsection*{Acknowledgment}
T.G. gratefully acknowledges the contribution of V. Oliker to unpublished collaborative work on an algorithm similar to the one described in this paper for a different reflector design problem. This work was done during 2002 to 2004 at Emory University.
\bibliographystyle{plain}
\bibliography{bibfile}

\end{document}